\documentclass[10pt, conference]{IEEEtran}
\IEEEoverridecommandlockouts
\usepackage{geometry}
\usepackage{cite}
\usepackage{amsmath,amssymb,amsfonts}
\usepackage{algorithmic}
\usepackage{graphicx}
\usepackage{textcomp}
\usepackage{xcolor}
\def\BibTeX{{\rm B\kern-.05em{\sc i\kern-.025em b}\kern-.08em
    T\kern-.1667em\lower.7ex\hbox{E}\kern-.125emX}}
\usepackage{subfigure}
\usepackage{subfigmat}
\usepackage{soul}
\usepackage{epsfig}

\usepackage{lipsum}

\usepackage{enumitem}

\graphicspath{{Picture/}}
\begin{document}

\title{\LARGE \bf  {Time-Fuel-Optimal Navigation of a Commercial Aircraft in Cruise with Heading and Throttle Controls using Pontryagin's Maximum Principle}}
%Characterization of the Optimal Controls in a Generic Cruise Phase w.r.t. Commercial Aircraft Flights 
%The Optimal Navigation in a Generic Cruise Phase w.r.t. Commercial Aircraft Flights
%Optimal Controls in Commercial Aircraft Flights for a Generic Cruise Model with Arbitrary Wind Fields
\author{Amin Jafarimoghaddam\thanks{Amin Jafarimoghaddam is with the Department of Aerospace Engineering, Universidad Carlos III de Madrid, Leganes, Madrid, Spain
		{\tt\small ajafarim@pa.uc3m.es}} and Manuel Soler\thanks{Manuel Soler is with the Department of Aerospace Engineering, Universidad Carlos III de Madrid, Leganes, Madrid, Spain
        {\tt\small masolera@ing.uc3m.es}}
}

\maketitle

\begin{abstract}
\textcolor{red}{In this research, we consider the commercial aircraft trajectory optimization problem for a general cruise model with arbitrary spatial wind fields to be solved using the Pontryagin's maximum principle. The model features two fundamental controls, namely "throttle setting" (which appears as a singular control) and "heading angle" (appearing as a regular control)}. For a constrained problem with minimum time-fuel objective, we show that the optimal heading angle is fully defined through the classic Zermelo's navigation identity. We also show that the optimal throttle setting can be characterized through a complete feedback function. The switching-point algorithm is employed to solve a case study where we inspect the optimality conditions and graph the optimal controls together with the optimal state and co-state variables.
\end{abstract}
\section{Introduction}
The agencies associated with commercial aircraft flights are inquisitively in demand of new optimization tools to minimize the direct operating cost, i.e., a combination of the fuel burn and the arrival time. \textcolor{red}{The optimization of commercial aircraft trajectories has been the subject of extensive research, with a variety of optimization techniques having been applied, however, the Pontryagin's maximum principle, has received comparatively little attention in the literature} \cite{1}. 

Here, we rely on the Pontryagin's maximum principle to solve minimum time-fuel commercial aircraft trajectory optimization problem for a general cruise model in the presence of an arbitrary spatial wind field.      

Focusing only on those approaches using the Pontryagin's maximum principle, the \textcolor{red}{commercial} aircraft navigation problem in the cruise phase with various cost functions has been solved in \cite{2, 3, 4, 5, 6, 7, 8}. More specifically, the range-optimal problem has been solved in \cite{2}. The same problem was revisited w.r.t. the compressibility effects in \cite{3}. Also, the fuel-optimal cruise at constant altitude with fixed arrival time has been solved in \cite{4}. Likewise, the fuel-optimal problem in a vertical plane (including the cruise phase) has been solved for structured flight segments in \cite{5}. The fuel-optimal problem in the cruise phase with a one-dimensional uniform wind field and fixed arrival time has been solved in \cite{6}. Approaches considering climate impact, such as \cite{7},  also exist in the literature, in which the contrails avoidance problem in the cruise phase was solved. Nevertheless, the study in \cite{7} ignores the speed dynamic and the derivative of the Hamiltonian w.r.t. the aircraft's mass and obtains a more general navigation formula. The direct operating cost with multiple cruise altitudes in the presence of wind has been addressed in \cite{8}. It is noteworthy that the speed dynamic has been ignored in the latter.   

\textcolor{red}{To summarize, the reviewed literature} (\cite{2, 3, 4, 5, 6, 7, 8}) \textcolor{red}{mainly focuses on simplified cruise models with only one active control (either "throttle setting", or "heading angle"). However, solving the optimal control problem associated with a general cruise model using the Pontryagin's maximum principle is particularly challenging due to the involvement of the two controls. In this research, we tackle this challenge and solve the problem for the first time.}

In this research, we consider a generic (realistic) cruise model, comprising of two active controls. The controls are: 1) the heading angle (a regular control), and, 2) the throttle setting (a singular control). The objective is to minimize the direct operating cost. We show that the optimal heading angle is fully characterized through the classic Zermelo's navigation identity even if the problem is subject to some standard state-inequality constraints (the classic Zermelo's navigation identity is the time-optimal, constant-speed navigation problem between two points inside a fluid flow \cite{9}). \textcolor{red}{For the optimal throttle setting, we employ successive derivatives of the switching function and leverage the optimality information related to the heading angle to characterize this singular control as a complete feedback function. In this respect, we show (through analyzing the switching function), that the optimal throttle setting is dependent on the optimal heading angle.} 

The switching-point algorithm \cite{10} is employed to solve a case study where the state-inequality constraints are inactive. In short, the switching-point algorithm runs a nonlinear programming only over the switching times and (possibly) some other scalar unknowns with a given control feedback. The switching-point algorithm, as in \cite{10}, is for a singular control problem without state-inequality constraints. Moreover, in \cite{10}, the decision variables can be the switching times and values of co-states at the entry of a singular arc. It is noteworthy that the study in \cite{10} is an extension to the previous works in this discipline (see e.g., \cite{11} and \cite{12}).  
\section{Problem Statement}\label{sec1}
The point-mass dynamics are commonly used to generate aircraft trajectories \cite{13}, \cite{14}. For commercial aircraft flights in cruise phase, the point-mass equations in Cartesian framework can be approximated as \cite{15}:
\begin{equation}\label{eq1}
	\begin{split}
	    &\frac{dx}{dt}=v(t)\cos{(\chi(t))}+w_x(x(t),y(t),h)=:F_x,\\
		&\frac{dy}{dt}=v(t)\sin{(\chi(t))}+w_y(x(t),y(t),h)=:F_y,\\
		&\frac{dv}{dt}=\frac{\Pi(t) T_{max}(h)-D(m(t),v(t),h)}{m(t)}=:F_v,\\
		&\frac{dm}{dt}=-\Pi(t) C_s(v(t))T_{max}(h)=:F_m.\\
	\end{split}
\end{equation}

In Eq. (\ref{eq1}), $x$, and $y$ are geometric variables, i.e., the cruise flight occurs in a horizontal $x-y$ plane, $v$ is the aerodynamic speed, $m$ is the aircraft mass, and $t$ is time. In addition, $w_x$ and $w_y$ are $x$, and $y$ components of wind respectively ($w_x$ and $w_y$ are known geometrical functions), $h$ is a constant altitude where the cruise flight occurs, $T_{max}$ is the \textcolor{red}{maximum} thrust force, $D$ is the drag force, and $C_s$ is the fuel flow. The controls are the heading angle ($\chi(t)$) and the throttle setting ($\Pi(t)$).   

%We simulate $T$, $D$, and $C_s$ by Base of Aircraft DAta (BADA3) model from EUROCONTROL providing a general smooth aircraft performance model []:
%\begin{equation}\label{eq2}
%	\begin{split}
%		&T(h)=C_{T_1}\big(1-\frac{h}{C_{T_2}}+h^2C_{T_3}\big),\quad C_{s}(v)=C_{s_1}\big(1+\frac{v}{C_{s_2}}\big)\\
%		&P(h)=P_0\big(\frac{\Theta_0-\beta h}{\Theta_0}\big)^{\frac{g}{\beta R}},\quad \rho(h)=\frac{P(h)}{R(\Theta_0-\beta h)}\\
%		&D(m,v,h)=\frac{1}{2}\rho(h)sv^2\big(C_{D_1}+C_{D_2}C_l^2\big)
%	\end{split}
%\end{equation}
%
%In Eq. (\ref{eq2}), $s$, $\Theta_0$, $\beta$, $R$, $g$, $P_0$, $C_{T_i}$, $i=1,2,3$, and $C_{D_i}$, $i=1,2$ are constants given in Table 1. Moreover, $C_l=\frac{2mg}{\rho(h)sv^2}$. 

The objective is:
\begin{equation}\label{eq2}
	\begin{split}
		\min_{\chi(t),\Pi(t),t_f}\mathcal{J}=\alpha t_f+(\alpha-1)m(t_f),\quad 0\leq\alpha\leq1.
	\end{split}
\end{equation}

Where $t_f$ denotes the final time.

\textcolor{red}{The cruise flight envelope is defined by the dynamic constraints described in Eq.} (\ref{eq1}) \textcolor{red}{, along with the following set of standard state-inequality constraints}:
\begin{equation}\label{eq3}
	\begin{split}
		&M_{min}\leq M(h,v(t))\leq M_{max},\\
		&v_{CAS,min}\leq v_{CAS}(h,v(t))\leq v_{CAS,max},\\
		&\quad\quad\quad\quad\forall t\in[t_0,t_f].
	\end{split}
\end{equation}

\textcolor{red}{and boundary conditions:}
\begin{equation}\label{eq3a}
	\begin{split}
		&x(t_f)=x_f,\quad y(t_f)=y_f,\quad v(t_f)=v_f,\\
		&x(t_0)=x_0,\quad y(t_0)=y_0,\quad v(t_0)=v_0,\quad m(t_0)=m_0.
	\end{split}
\end{equation}

In above, $t_0$ is the initial time, $M$ is the $Mach$ number, and $v_{CAS}$ is the calibrated airspeed \cite{16}. Moreover, $x_f,y_f,v_f,x_0,y_0,v_0,$ and $m_0$ are known values assigned for the boundary conditions.

The controls are also constrained as:
\begin{equation}\label{eq4}
	\begin{split}
		&\Pi_{min}\leq\Pi(t)\leq \Pi_{max},\\ &\chi_{min}\leq\chi(t)\leq\chi_{max},\\
		&\quad\quad \forall t\in[t_0,t_f].
	\end{split}
\end{equation}

It is noteworthy that, for the succeeding analysis, only the functionality of $w_x$, $w_y$, $T_{max}$, $D$, $C_s$, $Mach$, and $V_{CAS}$ is of relevance. \textcolor{red}{Nonetheless}, these terms will be elaborated in our case study (see sec.\ref{sec4}). 

\subsection{Compact Form Notation}\label{subsec2-1}
%Let us change the trigonometric terms in $F_x$, and $F_y$, by defining $q=tan(\chi(t))$. This gives $cos(\chi(t))=\frac{1}{\sqrt(1+q^2)}$, and $sin(\chi(t))=\frac{q}{\sqrt(1+q^2)}$.
The optimal control problem considered, can be written in a compact form notation as:
\begin{equation}\label{eq5a}
	\begin{split}
		&\min_{U(t),t_f} \mathcal{J}=\Phi({X}_f,t_f),\\
		&s.t.,\\
		&\frac{d{X}}{dt}=F(X(t),U(t)),\\
		&\mathcal{C}(U(t))\leq0,\quad \mathcal{S}(X(t))\leq0,\quad \forall t\in[t_0,t_f],\\
		&\phi_0(X_0)=0,\quad \phi_f(X_f)=0.\\
	\end{split}
\end{equation}

Where $X^T(t)=[x(t),y(t),v(t),m(t)]$, $F^T=[F_x,F_y,F_v,F_m]$, $\phi_0(X_0)=X_0-X(t_0)$, $\phi_0\in\mathbb{R}^4$, $\phi_f(X_f)=X_f-X(t_f)$, $\phi_f\in\mathbb{R}^3$, $U^T(t)=[\chi(t),\Pi(t)]$, and $\Phi(X_f,t_f)=\alpha t_f+(\alpha-1)m(t_f)$. The inequality constraints are:
\begin{equation}\label{eq6a}
	\begin{split}
		&\mathcal{C}(U(t))=
		\begin{pmatrix}
			\chi(t)-\chi_{max}\\
			\chi_{min}-\chi(t)\\
			\Pi(t)-\Pi_{max}\\
			\Pi_{min}-\Pi(t)
		\end{pmatrix}
	\end{split}.
\end{equation}
\begin{equation}\label{eq6b}
	\begin{split}
		&\mathcal{S}(X(t))=
		\begin{pmatrix}
			M(t)-M_{max}\\
			M_{min}-M(t)\\
			v_{CAS}(t)-v_{CAS,max}\\
			v_{CAS,min}-v_{CAS}(t)\\
		\end{pmatrix}.
	\end{split}
\end{equation}

By directly adjoining the constraints, we can define an augmented cost as (see e.g. \cite{17}):
\begin{equation}\label{eq5}
		\begin{split}
			&\bar{\mathcal{J}}:=\Phi(X_f,t_f)+\langle\nu_0,\phi_0(X_0)\rangle+\langle{\nu}_f,{\phi}_f(X_f)\rangle+\\&
			\int_{t_0}^{t_f}\bigg(\langle\lambda(t),F(X(t),U(t))-\frac{dX}{dt}\rangle+\langle\mu(t),\mathcal{C}(U(t))\rangle+\\
			&\langle\eta(t),\mathcal{S}(X(t))\rangle\bigg)dt.
		\end{split}
\end{equation}

In the above equation, the scalar multipliers are denoted by $\nu_0\in\mathbb{R}^4$, and $\nu_f\in\mathbb{R}^3$, while the inequality multipliers are represented by $\mu^T(t)=[\mu^\chi_u(t),\mu^\chi_l(t),\mu^{\Pi}_u(t),\mu^{\Pi}_l(t)]$, and $\eta(t)\in\mathbb{R}^4$. Additionally, the co-states are given by  $\lambda^T(t)=[\lambda_x(t),\lambda_y(t),\lambda_v(t),\lambda_m(t)]$. 

For this optimal control problem, the Hamiltonian is defined as \cite{17}:
\begin{equation}\label{eq7}
	\begin{split}
		&\mathcal{H}(X(t),U(t),\lambda(t),\mu(t),\eta(t)):=\\
		&\langle\lambda(t),F(X(t),U(t))\rangle+\langle\mu(t),\mathcal{C}(U(t))\rangle+\\
		&\langle\eta(t),\mathcal{S}(X(t))\rangle.
	\end{split}
\end{equation}

The optimality conditions are:
\begin{equation}\label{eq8}
	\begin{split}
		&\frac{\partial\mathcal{H}}{\partial U}=0,\\
		&\langle\mu(t),\mathcal{C}(U(t))\rangle=0,\quad \mu(t)\geq0,\quad \forall t\in[t_0,t_f],\\
	    &\langle\eta(t),\mathcal{S}(X(t))\rangle=0,\quad \eta(t)\geq0,\quad \forall t\in[t_0,t_f].\\
	\end{split}
\end{equation}

The co-state dynamics and the transversality conditions read:
\begin{equation}\label{eq9}
	\begin{split}
		&\frac{d{\lambda}^T}{dt}=-\frac{\partial\mathcal{H}}{\partial {X}},\\
		&{\lambda}_0:={\lambda}(t_0)=-\big[\frac{\partial{\phi}_0}{\partial {X}_0}\big]^T{\nu}_0,\\
		&{\lambda}_f:={\lambda}(t_f)=\frac{\partial\Phi}{\partial X_f}+\big[\frac{\partial{\phi}_f}{\partial {X}_f}\big]^T{\nu}_f,\\
		&\mathcal{H}_{f}:=\mathcal{H}(t_f)=-\frac{\partial\Phi}{\partial t_f}-{\nu}_f^T\frac{\partial{\phi}_f}{\partial t_f}.
	\end{split}
\end{equation}

Let $\tau$ be a possible time instant within the boundary arc at which the co-state variables are discontinuous. The jump conditions at the junction times read \cite{18}:
\begin{equation}\label{eq10}
	\begin{split}
		&\lambda^T(\tau^{-})=\lambda^T(\tau^{+})-\nu^T(\tau)\frac{\partial \mathcal{S}}{\partial X}|_{t=\tau},\\
		&\nu(\tau)\geq0,\quad \langle\nu(\tau),\mathcal{S}(X(\tau))\rangle=0.
	\end{split}
\end{equation}

\section{Classification of the Controls}\label{sec3}
The control problem defined through Eq. (\ref{eq1}) to Eq. (\ref{eq4}) is regular on the heading angle  ($\chi(t)$), and singular on the throttle setting ($\Pi(t)$). 

In order to classify the controls, we write:
\begin{equation}\label{eq11}
	\begin{split}
		&F(X(t),\chi(t),\Pi(t))=Q(X(t),\chi(t))+\Pi(t) P(X(t)).
	\end{split}
\end{equation}

In above:
\begin{equation}\label{eq12}
	\begin{split}
		&Q(X(t),\chi(t))=
		\begin{pmatrix}
			v(t)\cos(\chi(t))+w_x(x(t),y(t),h)\\
			v(t)\sin(\chi(t))+w_y(x(t),y(t),h)\\
			-\frac{D(m(t),v(t),h)}{m(t)}\\
			0
		\end{pmatrix},
	\end{split}
\end{equation}

and,
\begin{equation}\label{eq13}
	\begin{split}
		&P(X(t))=
		\begin{pmatrix}
			0\\
			0\\
			\frac{T_{max}(h)}{m(t)}\\
			-C_s(v(t))T_{max}(h)
		\end{pmatrix}.
	\end{split}
\end{equation}

\subsection{The Optimal $\chi(t)$ (Heading Angle)}\label{subsec31}
\textcolor{red}{We initially note that the state-inequality constraints} ($\mathcal{S}(X(t))$) \textcolor{red}{are functions of} $v$, \textcolor{red}{and} $h$. Since $h$ is a constant, from Eq. (\ref{eq10}), it is straightforward to show that the jump condition (upon existence), applies only to $\lambda_v(t)$.

The optimal $\chi(t)$, associated with the interior arc ( $\mu^\chi_l(t)=\mu^\chi_u(t)=0$), can be formulated through successive time-derivatives of the \textcolor{red}{first-order} optimality condition for $\chi(t)$. To this end, we write:  
%\begin{equation}\label{eq14}
%	\begin{split}
%    \frac{\partial\mathcal{H}}{\partial \chi}=0\rightarrow\frac{d^k}{dt^k}\frac{\partial\mathcal{H}}{\partial \chi}=0,\quad k=0,1,....
%	\end{split}
%\end{equation}
%
%We can check that Eq. (\ref{eq14}) for any $k$ involves only $\lambda_x$, and $\lambda_y$. \hl{Therefore by computing up to} $k=1$, \hl{we can extract a formula that only depends on the state variables.}

%In practice, $K$ is chosen in such a way to reach sufficient number of equations to remove the co-state variables. From Eq. (\ref{eq12}), and Eq. (\ref{eq13}), the variables $\chi(t)$, $x$, and $y$ appear only in $F_x$, and $F_y$. This means that Eq. (\ref{eq14}) for any $K$ only involves $\lambda_x$, and $\lambda_y$. Therefore, we can choose $K=1$, leading to a system of linear equations to remove the co-states $\lambda_x$, and $\lambda_y$.

%From Eq. (\ref{eq14}), for $k=0$ we can obtain:
\begin{equation}\label{eq15}
	\begin{split}
		&\frac{\partial\mathcal{H}}{\partial \chi}=0\rightarrow\langle\lambda(t),\frac{\partial Q}{\partial \chi}\rangle=0\rightarrow \tan(\chi(t))=\frac{\lambda_y(t)}{\lambda_x(t)}.\\
	\end{split}
\end{equation}

The first time-derivative of Eq. (\ref{eq15}) reads:
%\begin{equation}\label{eq16}
%	\begin{split}
%		&\langle\frac{d\lambda}{dt},\frac{\partial Q}{\partial q}\rangle+\langle\lambda,\frac{\partial}{\partial X}\frac{\partial Q}{\partial q}\big(Q+\Pi P(X)\big)\rangle+\\&\langle\lambda,\frac{\partial^2Q}{\partial q^2}\rangle\frac{dq}{dt}=0
%	\end{split}
%\end{equation} 
\begin{equation}\label{eq16}
	\begin{split}
		&\frac{d}{dt}\tan(\chi(t))=\frac{\frac{d\lambda_y}{dt}\lambda_x(t)-\frac{d\lambda_x}{dt}\lambda_y(t)}{\lambda_x^2(t)}.\\
	\end{split}
\end{equation}

From Eq. (\ref{eq9}),\textcolor{red}{the co-state dynamics} $\frac{d\lambda_y}{dt}$, and $\frac{d\lambda_x}{dt}$ \textcolor{red}{can be written as}:
%\begin{equation}\label{eq17}
%	\begin{split}
%		&\langle\lambda,\big[\frac{\partial Q}{\partial q},Q\big]\rangle+\langle\lambda,\big[\frac{\partial Q}{\partial q},P\big]\rangle\Pi(t)+\langle\lambda,\frac{\partial^2Q}{\partial q^2}\rangle\frac{dq}{dt}=0
%	\end{split}
%\end{equation}
\begin{equation}\label{eq17}
	\begin{split}
		&\frac{d\lambda_x}{dt}=-\big(\lambda_x(t)\frac{\partial w_x}{\partial x}+\lambda_y(t)\frac{\partial w_y}{\partial x}\big),\\
		&\frac{d\lambda_y}{dt}=-\big(\lambda_x(t)\frac{\partial w_x}{\partial y}+\lambda_y(t)\frac{\partial w_y}{\partial y}\big).
	\end{split}
\end{equation}

Plugging Eq. (\ref{eq17}) into Eq. (\ref{eq16}), and with the help of Eq. (\ref{eq15}), we obtain: 

%As stated earlier, Eq. (\ref{eq15}), and Eq. (\ref{eq17}), only involve $\lambda_x$, and $\lambda_y$. Moreover, from the definitions of $Q$, and $P$ (see Eq. (\ref{eq12}), and Eq. (\ref{eq13})), and by incorporating Eq. (\ref{eq15}), the middle term in Eq. (\ref{eq17}) vanishes. 
%
%Eq. (\ref{eq15}), and Eq. (\ref{eq17}) give:
%\begin{equation}\label{eq18}
%	\begin{split}
%		&\frac{dq}{dt}=\frac{A_1+qA_2}{B_1+qB_2}\\
%		&A=\big[\frac{\partial Q}{\partial q},Q\big]\quad B=\frac{\partial^2Q}{\partial q^2}	    
%	\end{split}
%\end{equation}
%
%In Eq. (\ref{eq18}), the subscripts stand for the vector elements.

%Eq. (\ref{eq18}), can be expanded and simplified as:
\begin{equation}\label{eq18}
	\begin{split}
		&\frac{d\chi}{dt}\big(1+\tan(\chi(t))^2\big)=-\frac{\partial w_x}{\partial y}+\\&\big(\frac{\partial w_x}{\partial x}-\frac{\partial w_y}{\partial y}\big)\tan(\chi(t))+\big(\frac{\partial w_y}{\partial x}\big)\tan(\chi(t))^2.	    
	\end{split}
\end{equation}

We observe that this equation is the Zermelo's navigation identity \cite{9}.

Eq. (\ref{eq18}) defines the optimal $\chi(t)$ for the interior arc even with active state-inequality constraints. Moreover, with the assumption that $\chi(t)$ is continuous in time, it will be straightforward to compute the boundary arc too.  

\subsection{The Optimal $\Pi(t)$ (Throttle Setting)}\label{subsec32}
Let us assume that the state-inequality constraints are inactive, i.e., $\eta(t)=0, \forall t\in[t_0,t_f]$.

The \textcolor{red}{first-order} optimality condition for $\Pi(t)$ reads:
\begin{equation}\label{eq19}
	\begin{split}
		&\frac{\partial\mathcal{H}}{\partial \Pi}=0\rightarrow \langle\lambda(t),P(X(t))\rangle+\mu^{\Pi}_u(t)-\mu^{\Pi}_l(t)=0.   
	\end{split}
\end{equation} 

we define the switching function as:
\begin{equation}\label{eq20}
	\begin{split}
		&S(t)=\langle\lambda(t),P(X(t))\rangle.  
	\end{split}
\end{equation} 

The bang-singular classification for the optimal $\Pi(t)$ is:
\begin{equation}\label{eq21}
	\begin{split}
		\Pi(t)=   
		\begin{cases}
			\Pi_{min} & S(t)>0,\\
			\Pi_{max} & S(t)<0,\\
			$undetermined$  & S(t)=0.
		\end{cases}
	\end{split}
\end{equation}

The optimal $\Pi(t)$ is called bang-bang if $S(t)$ has isolated zeros on an interval $I\subset[t_0,t_f]$; whereas the optimal $\Pi(t)$ is called singular if $S(t)=0$ holds for all $t\in I$.

For the singular $\Pi(t)$, since $S(t)=0$, we have: $S^{(k)}(t):=\frac{d^k}{dt^k}S(t)=0, k=1,2,...$. 

It can be checked that the singular $\Pi(t)$ appears in $S^{(2)}(t)$; that is to say, order of the singular arc is one. 

Since the Hamiltonian is not an explicit function of time, we have: $\mathcal{H}(t)=constant$ ($\forall t\in[t_0,t_f]$). Therefore, from Eq. (\ref{eq9}), we have: $\mathcal{H}(t)=-\alpha$. Moreover, $\frac{\partial\mathcal{H}}{\partial \chi}=0$ holds at all times. 

On using Eq. (\ref{eq11}), the co-state dynamics can be written as:
\begin{equation}\label{eq221}
	\begin{split}
		&\frac{d\lambda^T}{dt}=-\lambda^T(t)\big(\frac{\partial Q}{\partial X}+\Pi(t)\frac{\partial P}{\partial X}\big).
	\end{split}
\end{equation}

The first time-derivative of the switching function is computed as:
\begin{equation}\label{eq222}
	\begin{split}
		&S^{(1)}(t)=\langle\frac{d\lambda}{dt},P(X(t))\rangle+\langle\lambda(t),\frac{dP}{dt}\rangle=0.
	\end{split}
\end{equation}

Noting that $\frac{dP}{dt}=\frac{\partial P}{\partial X}\frac{dX}{dt}$, it is straightforward to show that $\Pi(t)$ drops from $S^{(1)}(t)$. As a result, $S^{(1)}(t)$ can be expressed as follows:
\begin{equation}\label{eq223}
	\begin{split}
		&S^{(1)}(t)=\langle\lambda(t),\mathcal{A}(X(t),\chi(t))\rangle.
	\end{split}
\end{equation}

In Eq. (\ref{eq223}), the vector $\mathcal{A}$ is:
\begin{equation}\label{eq224}
	\begin{split}
		&\mathcal{A}(X(t),\chi(t))=\frac{\partial P}{\partial X}Q(X(t),\chi(t))-\frac{\partial Q}{\partial X}P(X(t)).
	\end{split}
\end{equation} 

\textcolor{red}{It should be noted that the vector} $\mathcal{A}$ \textcolor{red}{can also be represented using Lie bracket notations}.

%On using the co-state dynamics $\frac{d\lambda^T}{dt}=-\lambda^T(t)\big(\frac{\partial Q}{\partial X}+\Pi(t)\frac{\partial P}{\partial X}\big)$, the set of algebraic equations defining the co-states can be written as:
\textcolor{red}{We exploit the following set of algebraic equations to express the co-states in explicit terms:}
\begin{equation}\label{eq22}
	\begin{split}
		&S(t)=0\rightarrow\langle\lambda(t),P(X(t))\rangle=0,\\
		&S^{(1)}(t)=0\rightarrow\langle\lambda(t),\mathcal{A}(X(t),\chi(t))\rangle=0,\\
		&\mathcal{H}(t)=-\alpha\rightarrow\langle\lambda(t),Q(X(t),\chi(t))\rangle=-\alpha,\\
		&\frac{\partial\mathcal{H}}{\partial \chi}=0\rightarrow\langle\lambda(t),\frac{\partial Q}{\partial \chi}\rangle=0.
	\end{split}
\end{equation}

%In Eq. (\ref{eq22}), the vector $\mathcal{A}$ is:
%\begin{equation}\label{eq22a}
%	\begin{split}
%		&\mathcal{A}(X(t),\chi(t))=\frac{\partial P}{\partial X}Q(X(t),\chi(t))-\frac{\partial Q}{\partial X}P(X(t)).
%	\end{split}
%\end{equation} 
%
%\hl{It should be noted that the vector} $\mathcal{A}$ \hl{can also be represented using Lie bracket notations}.

The above algebraic linear system can be simply solved by any platform supporting symbolic computations such as MATLAB. To this end, we write the solution to system (\ref{eq22}) as:
\begin{equation}\label{eq23}
	\begin{split}
		\det[\bar{\bar{\mathcal{M}}}]	
		\begin{pmatrix}
			\lambda_x(t)\\
			\lambda_y(t)\\
			\lambda_v(t)\\
			\lambda_m(t)
		\end{pmatrix}
		=
		adj[\bar{\bar{\mathcal{M}}}]
		\mathcal{R}.
	\end{split}
\end{equation}

Where:
\begin{equation}\label{eq24}
	\begin{split}	
		\bar{\bar{\mathcal{M}}}
		=
		\begin{pmatrix}
			0&0&P_3&P_4\\
			\mathcal{A}_1&\mathcal{A}_2&\mathcal{A}_3&\mathcal{A}_4\\
			Q_1&Q_2&Q_3&0\\
			\tan(\chi(t))&-1&0&0
		\end{pmatrix},
	\end{split}
\end{equation}

and,
\begin{equation}\label{eq25}
	\begin{split}	
		\mathcal{R}
		=
		\begin{pmatrix}
			0\\
			0\\
			-\alpha\\
			0
		\end{pmatrix}.
	\end{split}
\end{equation}

In Eq. (\ref{eq24}), subscripts stand for the vector elements.

In order to obtain the optimal singular $\Pi(t)$, we write:
\begin{equation}\label{eq26a}
	\begin{split}
		&S^{(2)}(t)=\frac{d}{dt}S^{(1)}(t)=\frac{d}{dt}\langle\lambda(t),\mathcal{A}(X(t),\chi(t))\rangle=\\
		&\langle\frac{d\lambda}{dt},\mathcal{A}(X(t),\chi(t)\rangle+\langle\lambda(t),\frac{d\mathcal{A}}{dt}\rangle=0.
	\end{split}
\end{equation}

\textcolor{red}{The time-derivative of} $\mathcal{A}$ \textcolor{red}{is computed as}:
\begin{equation}\label{eq26b}
	\begin{split}
		&\frac{d\mathcal{A}}{dt}=\frac{\partial\mathcal{A}}{\partial X}\frac{dX}{dt}+\frac{\partial\mathcal{A}}{\partial \chi}\frac{d\chi}{dt}.
	\end{split}
\end{equation}

With the help of Eq. (\ref{eq26a}), and Eq. (\ref{eq26b}), together with the co-state dynamics (Eq. (\ref{eq221})), one arrives at:
\begin{equation}\label{eq26}
	\begin{split}
		&S^{(2)}(t)=0\rightarrow\Pi(t)=-\frac{\langle\lambda(t),\mathcal{B}\rangle+\langle\lambda(t),\frac{\partial\mathcal{A}}{\partial \chi}\rangle\frac{d\chi}{dt}}{\langle\lambda(t),\mathcal{D}\rangle}.
	\end{split}
\end{equation}

In above, the vectors $\mathcal{B}$, and $\mathcal{D}$ are:
\begin{equation}\label{eq26c}
	\begin{split}
		&\mathcal{B}=\frac{\partial\mathcal{A}}{\partial X}Q(X(t),\chi(t))-\frac{\partial Q}{\partial X}\mathcal{A}(X(t),\chi(t)),\\
		&\mathcal{D}=\frac{\partial\mathcal{A}}{\partial X}P(X(t))-\frac{\partial P}{\partial X}\mathcal{A}(X(t),\chi(t)).\\
	\end{split}
\end{equation}

%\big[\big[P,Q\big],Q\big]
%\big[\big[P,Q\big],P\big]
The generalized Legendre-Clebsch (LC) second-order necessary conditions dictate \cite{19}:
\begin{equation}\label{eq27}
	\begin{split}
		&-\langle\lambda(t),\mathcal{D}\rangle\geq0,\quad \forall t\in\Omega_s.
	\end{split}
\end{equation}

Where $\Omega_s$ is an interval where the optimal $\Pi(t)$ is singular\footnote{\textcolor{red}{In the current study, all computations related to the optimal} $\chi(t)$ \textcolor{red}{and singular} $\Pi(t)$ \textcolor{red}{were carried out symbolically in MATLAB. This involved solving the co-state system defined by Eq.} (\ref{eq23}) \textcolor{red}{and deriving symbolic formulas for} $\mathcal{A}$, $\mathcal{B}$, and $\mathcal{D}$.}.

\subsubsection{The Special Case $\alpha=0$}\label{subsubsec321}
Since $\mathcal{R}=\vec{0}$ if $\alpha=0$, the solution to Eq. (\ref{eq23}) becomes intractable. On this occasion, from Eq. (\ref{eq23}), we have:
\begin{equation}\label{eq28}
	\begin{split}
		&\det\big[\bar{\bar{\mathcal{M}}}\big]=0\rightarrow\frac{d}{dt}\det\big[\bar{\bar{\mathcal{M}}}\big]=0.
	\end{split}
\end{equation}

Therefore, from Eq. (\ref{eq28}), we obtain the optimal singular $\Pi(t)$ as:
\begin{equation}\label{eq29}
	\begin{split}
		&\Pi(t)=-\frac{\frac{\partial}{\partial X}\det\big[\bar{\bar{\mathcal{M}}}\big]Q(X(t),\chi(t))}{\frac{\partial}{\partial X}\det\big[\bar{\bar{\mathcal{M}}}\big]P(X(t))}.
	\end{split}
\end{equation}

It is noteworthy that we can also handle the case $\alpha=0$ asymptotically, i.e., to compute for $\alpha\rightarrow 0$.

\section{Case Study}\label{sec4}
In accordance with our analysis presented in the previous sections, we consider inequality constraints only w.r.t. the throttle setting, i.e., $\Pi_{min}\leq\Pi(t)\leq\Pi_{max},\forall t\in[0,t_f]$.
 
We use the BADA3 model for $T_{max}(h)$, $C_s(v)$, and $D(m,v,h)$ \cite{20}. In addition, air density is approximated by International Standard Atmospheric (ISA) model: 
\begin{equation}\label{eq30}
\begin{split}
	&T_{max}(h)=C_{T_1}\big(1-\frac{h}{C_{T_2}}+h^2C_{T_3}\big),\\
	&P(h)=P_0\big(\frac{\Theta_0-\beta h}{\Theta_0}\big)^{\frac{g}{\beta R}},\quad \rho(h)=\frac{P(h)}{R(\Theta_0-\beta h)},\\
	&D(m,v,h)=\frac{1}{2}\rho(h)sv^2\big(C_{D_1}+C_{D_2}C_l^2\big),\\
	&C_l=\frac{2mg}{\rho s v^2},\quad C_{s}(v)=C_{s_1}\big(1+\frac{v}{C_{s_2}}\big).
\end{split}
\end{equation}

In above, $s$ is the aerodynamic lift surface and $\rho$ is the air density. $C_{T_i}, i=1,2,3$, $s$, $C_{D_i}, i=1,2$, $C_{s_i}, i=1,2$, $R$, $\beta$, $P_0$, $g$, and $\theta_0$ are known constants for a medium-haul aircraft (see \cite{16} for a quick access to their specific values and definitions). 

\textcolor{red}{Since} $h$ \textcolor{red}{is a constant (i.e., the altitude where the cruise flight occurs), the maximum thrust force} $T_{max}(h)$ \textcolor{red}{will be a constant too}.
%\begin{table}[hbt!]\label{table1}
%	\caption{Aircraft model constants}
%	\begin{center}\label{table1}
%		\begin{tabular}{|c||c||c|}
%			\hline 
%			Parameter & Value & Unit (SI)\\
%			\hline
%			$s$ &122.6&$m^2$\\
%			\hline
%			$C_{T_1}$&141040&$N$\\
%			\hline
%			$C_{T_2}$&14909.9&$m$\\
%			\hline
%			$C_{T_3}$&$6.997\times10^{-10}$&$m^{-2}$\\
%			\hline
%			$ C_{D_1}$&0.0242&$-$\\
%			\hline
%			$C_{D_2}$&0.0469&$-$\\
%			\hline
%			$C_{s_1}$&$1.055\times10^{-5}$&$kg.s^{-1}N^{-1}$\\
%			\hline
%			$C_{s_2}$&441.54&$m.s^{-1}$\\
%			\hline
%			$\Theta_0$&288.15&$K$\\
%			\hline
%			$P_0$&101325&$Pa$\\
%			\hline
%			$g$&9.81&$m/s^2$\\
%			\hline
%			$R$&287.058&$J.kg^{-1}K^{-1}$\\
%			\hline
%			$\beta$&0.0065&$K.m^{-1}$\\
%			\hline
%		\end{tabular}
%	\end{center}
%\end{table}   

The wind components are simulated by a general second-order polynomial for each component.

Without loss of generality, we can consider $w_x=w_x(x,y)$, and $w_y=w_y(x,y)$, and drop $h$. We also note that the vertical component of wind is assumed to be zero \cite{15}, \cite{13}. Therefore, for the sake of consistency with the fluid flow behavior, the continuity equation must be satisfied, i.e., $\frac{\partial w_x}{\partial x}+\frac{\partial w_y}{\partial y}=0$. With this observation, the wind components become:
\begin{equation}\label{eq31}
 	\begin{split}
 		&w_x(x,y)=(a_0\bar{w}_x^b)+(a_1\frac{\bar{w}_x^b}{x_f})x+(a_2\frac{\bar{w}_x^b}{x_f^2})x^2+\\&(a_3\frac{\bar{w}_x^b}{y_f})y+(a_4\frac{\bar{w}_x^b}{y_f^2})y^2+(a_5\frac{\bar{w}_x^b}{x_fy_f})xy,\\
 		&w_y(x,y)=-\int\frac{\partial w_x}{\partial x}dy+f(x)\\
 		&f(x)=\bar{w}_y^b+(b_0\frac{\bar{w}_y^b}{x_f})x+(b_1\frac{\bar{w}_y^b}{x_f^2})x^2.
 	\end{split}
\end{equation}

In above, $a_i, i=0,..,5$, and $b_i, i=0,1$ are in general dimensionless random values between $-1$ and $+1$. Moreover, $\bar{w}_x^b$, and $\bar{w}_y^b$ are average dimensional wind constants. 

\textcolor{red}{Upon conducting an estimation of real wind data (w.r.t. various, though small, atmospheric zones) using the aforementioned second-order model, we observed a very good level of consistency. In particular, the total modeling error for the wind components was found to be below 10 percent in relation to low-resolution wind data. It is worth noting that in this study, we adopt a generalized approach and do not fine-tune the model based on specific data sets.}
   
The tabulated parameters in table (\ref{table2}) are those we have fixed in our simulations. Therefore, one can check that the only free parameter in our simulations is $\alpha$.
\begin{table}[hbt!]\label{table2}
	\caption{boundary conditions, bounds, and the selected snapshot of the wind parameters}
	\begin{center}\label{table2}
		\begin{tabular}{|c||c||c||c|} 
%			\hline
%			Parameter & Value&Parameter & Value\\
			\hline
			$x_0$ &0$(m)$&$x_f$&1.5$\times 10^6(m)$\\
			\hline
			$y_0$&0$(m)$&$y_f$&7$\times 10^5(m)$\\
			\hline
			$v_0$&200$(m/s)$&$v_f$&200$(m/s)$\\
			\hline
			$m_0$&59000$(kg)$&$h$&10000$(m)$\\
			\hline
			$\Pi_{max}$&1&$\Pi_{min}$&0\\
			\hline
			$\bar{w}_x^b$&40$(m/s)$&$\bar{w}_y^b$&-20$(m/s)$\\
			\hline
			$a_0$&0.77406&$a_1$&-0.86240\\
			\hline
			$a_2$&-0.63294&$a_3$&0.47414\\
			\hline
			$a_4$&0.39342&$a_5$&0.55398\\
			\hline
			$b_0$&0.00380&$b_1$&-0.14900\\
			\hline
		\end{tabular}
	\end{center}
\end{table}
			
\section{Computational Algorithm}
Preliminary analysis (\textcolor{red}{using a single-shooting Euler-based direct transcription method with high number of grids}) implies that the optimal $\Pi(t)$ contains at most one (interior) singular arc between two boundary arcs. 

We have employed the switching-point algorithm, taking the switching times and the initial heading as decision variables (see \cite{10} for more information about the switching-point algorithm and the associated mathematical justifications).

More specifically, we extend the state dynamics by Eq. (\ref{eq18}), and the nonlinear programming becomes:
\begin{equation}\label{eq32}
	\begin{split}
		&\min_{\chi(0), t_1, t_2, t_f} \hat{\mathcal{J}}=\Phi(X_f,t_f),\\
		&s.t.,\\
		&\frac{d{X}}{dt}=F(X(t),\chi(t),\Pi(t)),\\
		&\frac{d\chi}{dt}\big(1+\tan(\chi(t))^2\big)=-\frac{\partial w_x}{\partial y}+\\&\big(\frac{\partial w_x}{\partial x}-\frac{\partial w_y}{\partial y}\big)\tan(\chi(t))+\big(\frac{\partial w_y}{\partial x}\big)\tan(\chi(t))^2,\\
		&{\phi}_0(X_0)=0,\\
		&{\phi}_f(X_f)=0,\\
		&0\leq t_1\leq t_2\leq t_f.
	\end{split}
\end{equation}

With:
\begin{equation}\label{eq33}
	\begin{split}
		\Pi(t)=   
		\begin{cases}
			\Pi_{max} & t\in[0,t_1),\\
			Eq. (\ref{eq26}) \lor Eq. (\ref{eq29})  & t\in[t_1,t_2],\\
			\Pi_{min} & t\in(t_2,t_f].	
		\end{cases}
	\end{split}
\end{equation}

We have adopted the interior-point/barrier algorithm of the nonlinear programming solver $fmincon$ from MATLAB® Optimization Toolbox as the optimization module.

We use the results due to the above nonlinear programming to compute the co-state variables (and accordingly, the switching function) over the boundary arcs. This is doable by backward and forward integration of the co-state dynamics from $t_1$, and $t_2$ respectively.

\subsection{Special Case with a Constant Wind Field}\label{subsecnew1} 
\textcolor{red}{Assuming a constant wind field, it is demonstrated that the optimization problem formulated in Eq.} (\ref{eq32}) \textcolor{red}{can exclude} $\chi(0)$ \textcolor{red}{as a decision variable.}

Suppose that: $w_x=constant=:W_x$, and $w_y=constant=:W_y$. 
%the following holds: $\tan(\chi_0)=\frac{y_f-W_yt_f-y_0}{x_f-W_xt_f-x_0}$.

From the system dynamics Eq. (\ref{eq1}), we can write:
\begin{equation}\label{eq34}
	\begin{split}
		&\frac{dx}{dt}=v(t)\cos(\chi(t))+W_x,\quad\frac{dy}{dt}=v(t)\sin(\chi(t))+W_y.
	\end{split}
\end{equation}

Since $W_x$, and $W_y$ are constants, from Eq. (\ref{eq18}) we have: $\chi(t)=\chi(0)=constant$. Therefore, by integrating Eq. (\ref{eq34}) from $0$ to $t_f$, and after some elementary manipulations, we get: 
\begin{equation}\label{eq35}
	\begin{split}
		&\tan(\chi(0))=\frac{y_f-W_yt_f-y_0}{x_f-W_xt_f-x_0}.
\end{split}
\end{equation}

\textcolor{red}{Therefore, in case of having a constant wind field, the initial heading angle} $\chi(0)$ \textcolor{red}{is a function of the constant wind components, boundary conditions, and the final time} $t_f$.
  
\section{Numerical Results}
\begin{figure}[hbt!]	
	{\includegraphics[width=1\columnwidth]{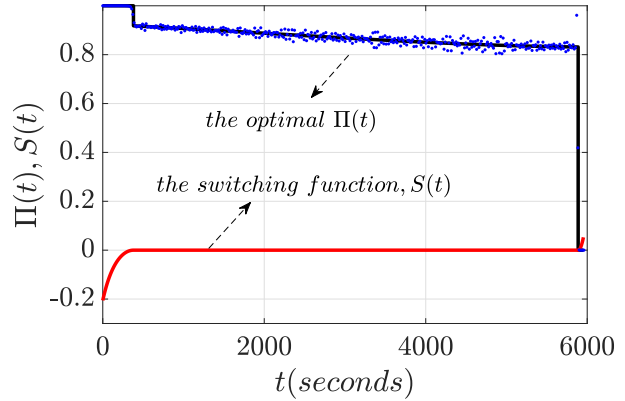}}   \hspace{-3cm}
	\caption{The optimal $\Pi(t)$ (black), and $S(t)$ (red) where $\alpha=0.4$, compared to single-shooting Euler-based direct transcription method with 400 nodes (blue dots); The optimal cost by the direct method:-30109.35, The optimal cost by the indirect method:-30109.38} 
	\label{Fig. 1}
\end{figure}

We have obtained optimal results for various values of $\alpha$. From the definition of the cost function, $\alpha$ determines the trade-off between fuel-optimal and time-optimal problems. For each studied $\alpha$, we have checked the second-order optimality condition (Eq. (\ref{eq27})) for the singular arc and the first-order optimality condition for the boundary arcs. In addition, for each stage of $\alpha$, we have compared our results with the results due to a \textcolor{red}{single-shooting Euler-based direct transcription method with high number of grids} (see Fig. \ref{Fig. 1}). \textcolor{red}{Upon examining this figure, we can ascertain that the optimal costs derived from the direct method and the present indirect method exhibit remarkable similarity. However, it is noteworthy that the estimated singular arc generated by the direct method exhibits chattering behavior around the singular solution by the indirect method.}

Moreover, we note that the switching-point algorithm does not directly account for $\lambda(t_f)$. From the transversality conditions, we can check that $\lambda_m(t_f)=\alpha-1$. By computing the co-state variables over the boundary arcs (after nonlinear programming), we have checked that $|\lambda_m(t_f)-(\alpha-1)|<10^{-4}$. This also stands as an additional confirmation of the obtained optimal results (see Fig. \ref{Fig. 8}).

\begin{figure}[hbt!]	
	{\includegraphics[width=1\columnwidth]{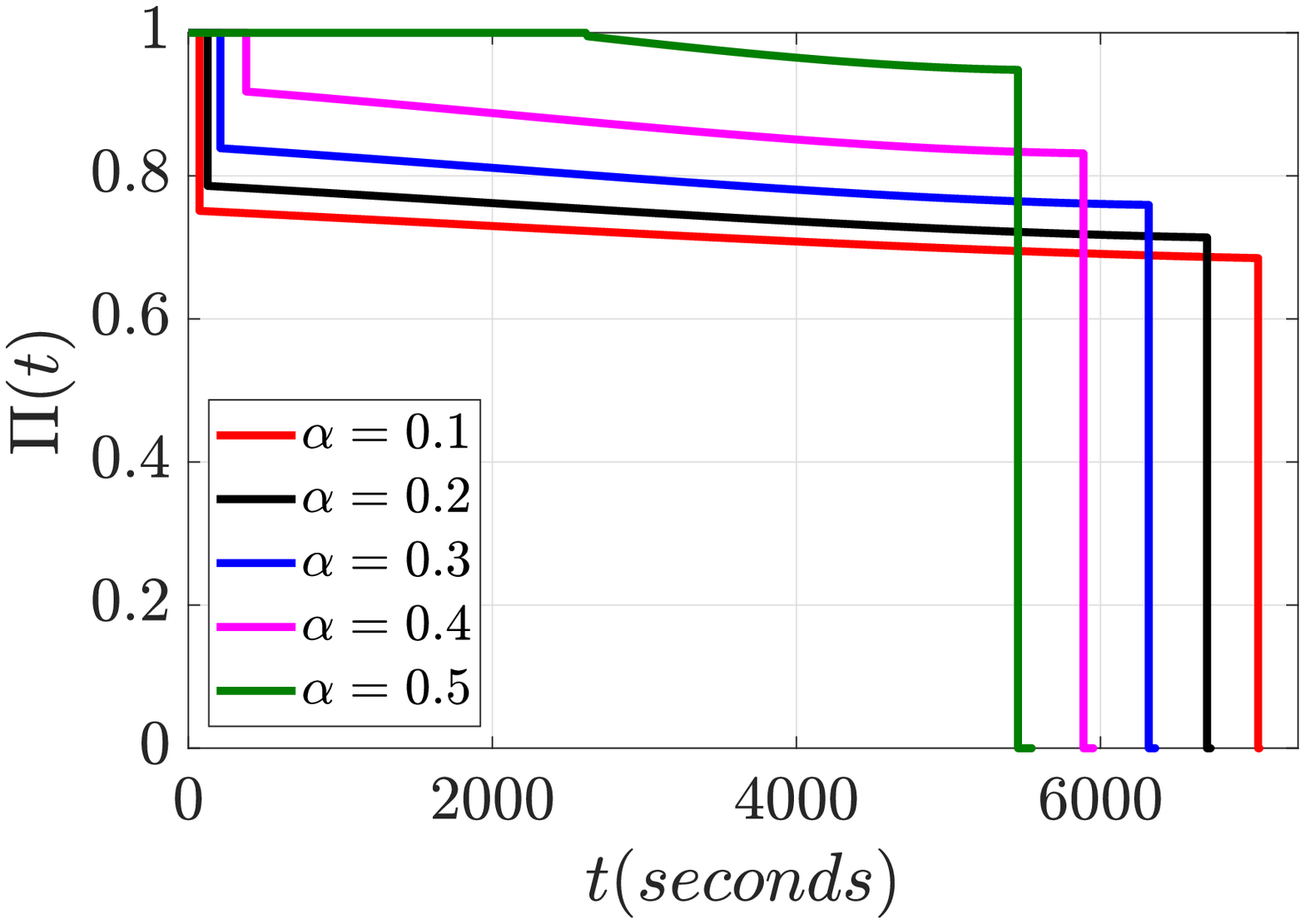}}   \hspace{-3cm}
	\caption{The optimal $\Pi(t)$ in different $\alpha$} 
	\label{Fig. 2}
\end{figure}

\begin{figure}[hbt!]	
	{\includegraphics[width=1\columnwidth]{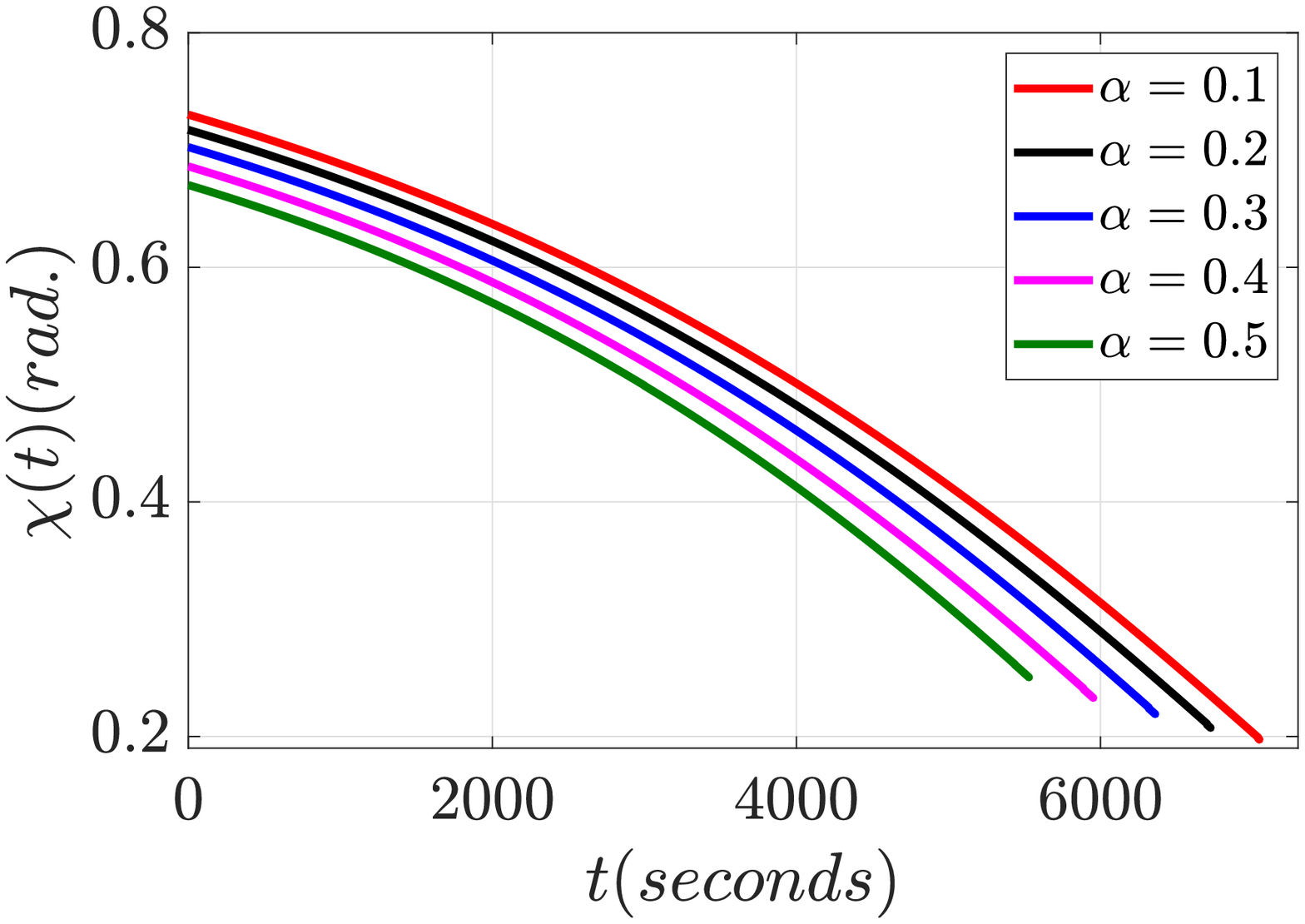}}   \hspace{-3cm}
	\caption{The optimal $\chi(t)$ in different $\alpha$} 
	\label{Fig. 3}
\end{figure}

\begin{figure}[hbt!]	
	{\includegraphics[width=1\columnwidth]{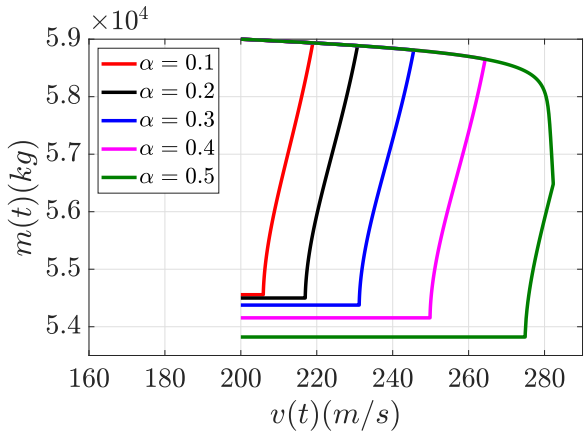}}   \hspace{-3cm}
	\caption{The optimal $m(t)$ as a function of $v(t)$ in different $\alpha$} 
	\label{Fig. 4}
\end{figure}

The optimal controls $\Pi(t)$, $\chi(t)$, and the optimal states ($m(t)$ as a function of $v(t)$) are shown in Fig. \ref{Fig. 2}, Fig. \ref{Fig. 3}, and Fig. \ref{Fig. 4} respectively. \textcolor{red}{Based on} Fig. \ref{Fig. 2} and Fig. \ref{Fig. 4}, \textcolor{red}{it is clear that an increase in} $\alpha$ \textcolor{red}{leads to a corresponding increase in the optimal speed by adjusting the throttle setting.} 

From Fig. \ref{Fig. 2}, \textcolor{red}{it can be observed that as} $\alpha$ \textcolor{red}{increases, the first "bang" segment of the optimal throttle} $\Pi(t)$ \textcolor{red}{expands. This implies that there is a specific value of} $\alpha$ at which the optimal $\Pi(t)$ \textcolor{red}{switches to a "bang-bang" control.}. 

\begin{figure}[hbt!]	
	{\includegraphics[width=1\columnwidth]{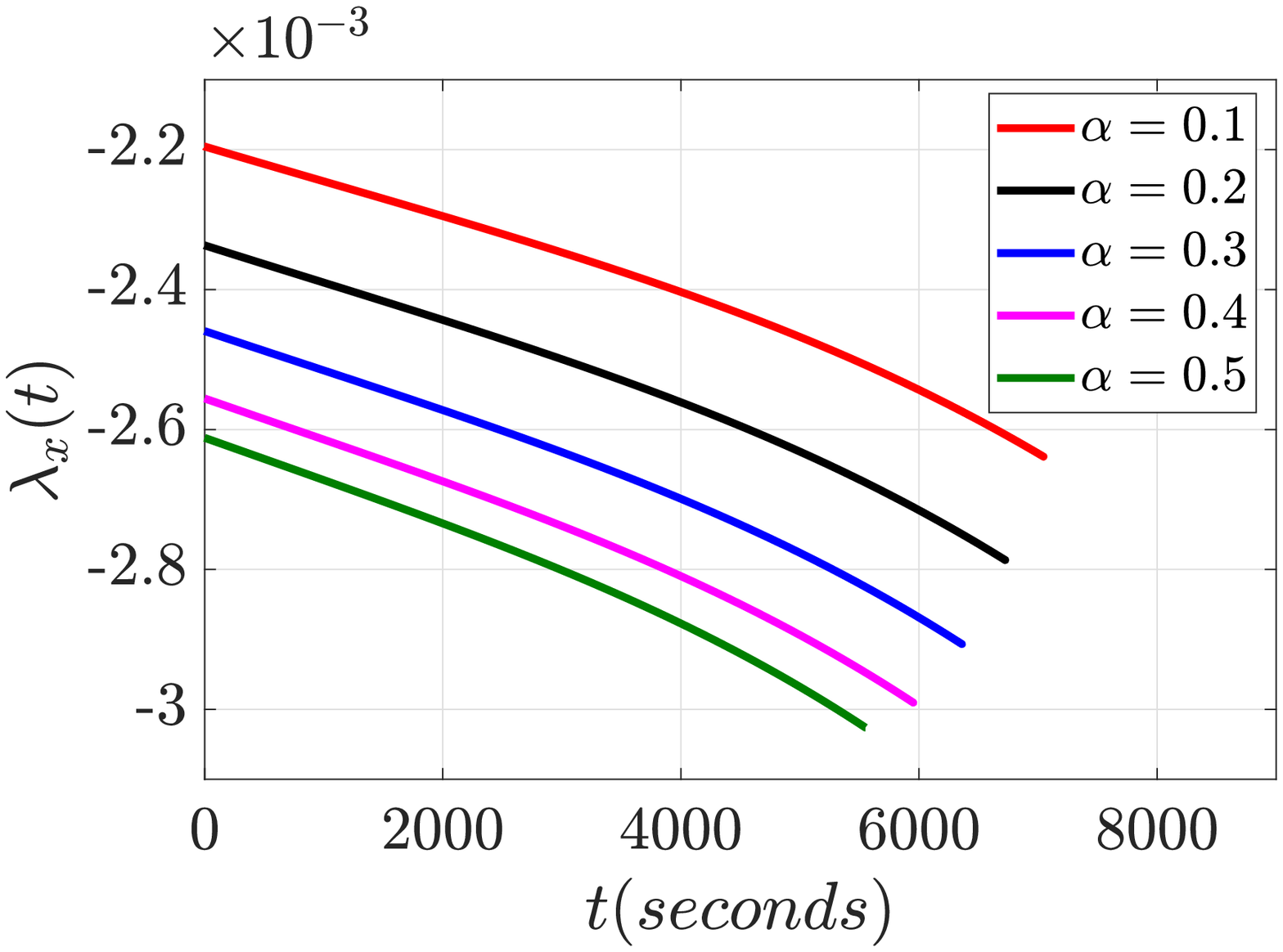}}   \hspace{-3cm}
	\caption{The optimal $\lambda_x$ in different $\alpha$} 
	\label{Fig. 5}
\end{figure}

\begin{figure}[hbt!]	
	{\includegraphics[width=1\columnwidth]{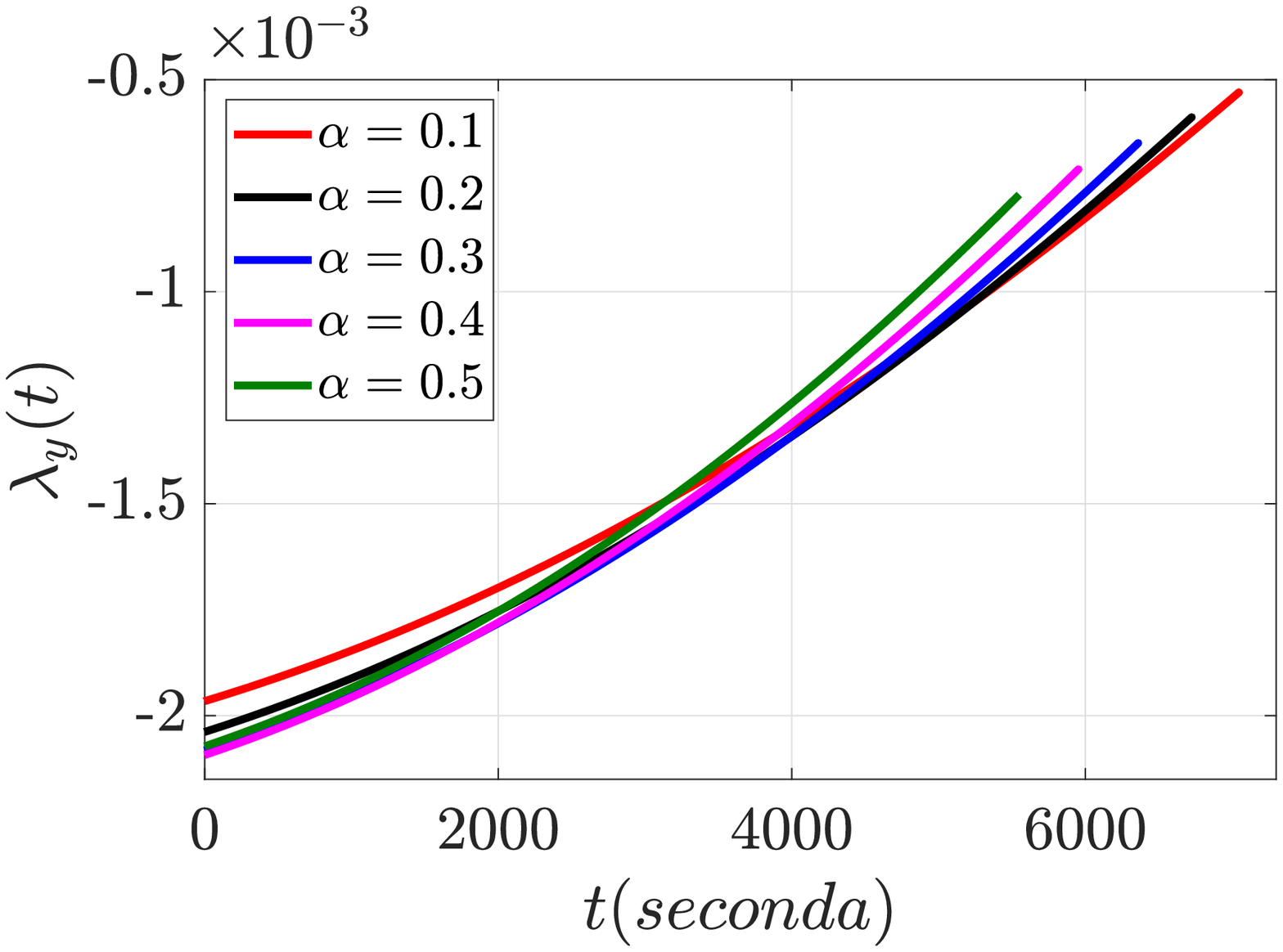}}   \hspace{-3cm}
	\caption{The optimal $\lambda_y$ in different $\alpha$} 
	\label{Fig. 6}
\end{figure}

\begin{figure}[hbt!]	
	{\includegraphics[width=1\columnwidth]{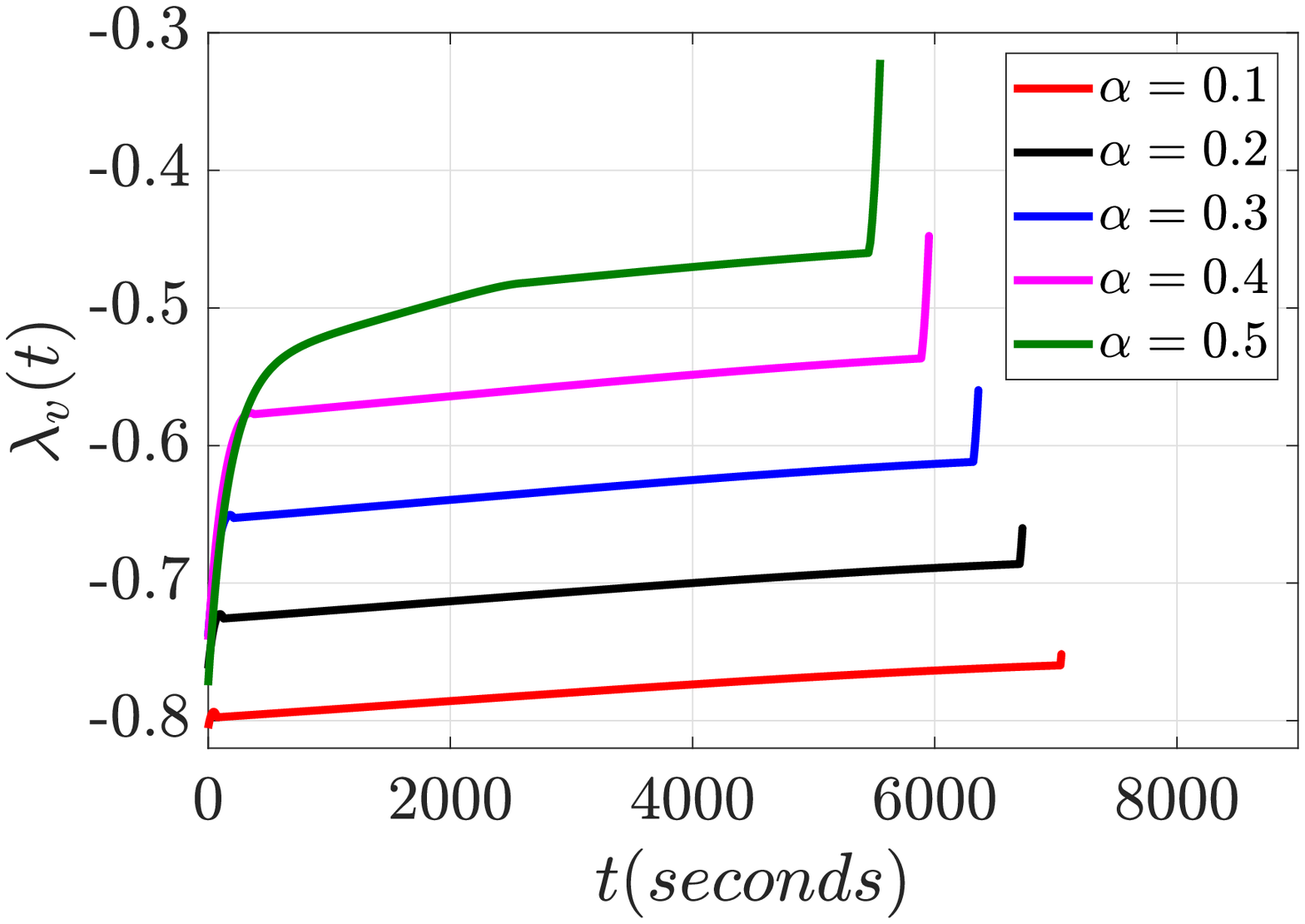}}   \hspace{-3cm}
	\caption{The optimal $\lambda_v$ in different $\alpha$} 
	\label{Fig. 7}
\end{figure}

\begin{figure}[hbt!]	
	{\includegraphics[width=1\columnwidth]{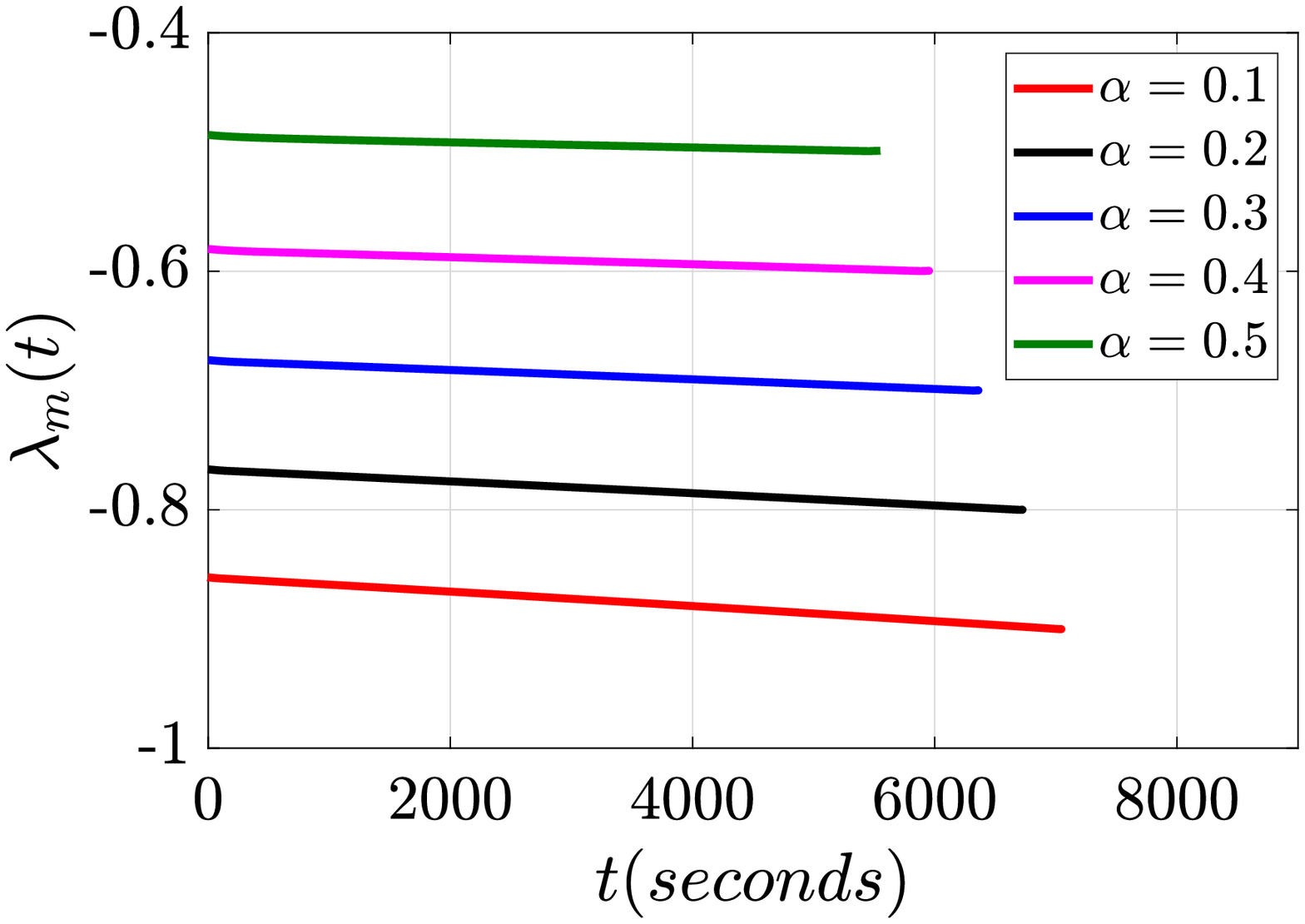}}   \hspace{-3cm}
	\caption{The optimal $\lambda_m$ in different $\alpha$} 
	\label{Fig. 8}
\end{figure}

The optimal co-state variables in different values of $\alpha$ are graphed in Fig. \ref{Fig. 5}-Fig. \ref{Fig. 8}.

%From Fig. \ref{Fig. 5}-Fig. \ref{Fig. 8}, it appears that the objective function exhibits more sensitivity to changes in aircraft's speed in comparison to the other state variables. 
As depicted in Fig. \ref{Fig. 8}, $\lambda_m(t)$ \textcolor{red}{displays only marginal changes with respect to time. Hence, it may be reasonable to consider it as a constant for analysis purposes, especially in engineering applications where rough estimations of optimality are sufficient.}

\begin{figure}[hbt!]	
	{\includegraphics[width=1\columnwidth]{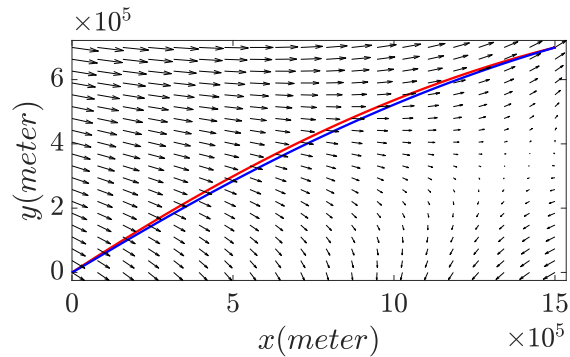}}   \hspace{-3cm}
	\caption{The optimal $x-y$ {trajectories (lateral path)} in different $\alpha$: blue and red lines show the lateral paths w.r.t. $\alpha=0.5$, and $\alpha=0.1$ respectively.} 
	\label{Fig. 9}
\end{figure}

Fig. \ref{Fig. 9} \textcolor{red}{illustrates the impact of} $\alpha$ \textcolor{red}{on the optimal $x-y$ trajectories}. \textcolor{red}{Upon examination of the figure, it is apparent that the parameter} $\alpha$ \textcolor{red}{exhibits negligible influence on the evolution of the} $x-y$ \textcolor{red}{trajectories. More precisely, the} $x-y$ \textcolor{red}{trajectories are more contingent on the wind configuration.}



\section{Conclusion and Future Works}
\textcolor{red}{Pontryagin's maximum principle was applied to solve a general (realistic) version of the optimization problems related to commercial aircraft trajectory in cruise phase. The analysis focused on the control functions, namely the "heading angle" (regular control) and the "throttle setting" (singular control), from which optimality formulas were derived. To handle the singular control, the switching-point algorithm was utilized as an alternative approach to the conventional shooting methods}. We have observed that the singular $\Pi(t)$ vanishes in larger values of $\alpha$. Therefore, an open line of research is to explore the condition in which the singular $\Pi(t)$ disappears. Moreover, our case study did not involve state-inequality constraints. This can also be a subject for the future research.

%%%%%%%%%%%%%%%%%%%%%%%%%%%%%%%%%%%%%%%%%%%%%%%%%%%%%%%%%%%%%%%%%%%%%%%%%%%%%%%
%\begin{thebibliography}{99}
%
%\bibitem{c1}
%J.G.F. Francis, The QR Transformation I, {\it Comput. J.}, vol. 4, 1961, pp 265-271.
%
%\bibitem{c2}
%Douglas M. Pargett and Mark D. Ardema. "Flight path optimization at constant altitude", {\it Journal of Guidance, Control, and Dynamics}, 30(4):1197-1201, 2007.
%
%\bibitem{c3}
%D. Boley and R. Maier, "A Parallel QR Algorithm for the Non-Symmetric Eigenvalue Algorithm", {\it in Third SIAM Conference on Applied Linear Algebra}, Madison, WI, 1988, pp. A20.
%
%\end{thebibliography}

%%%%%%%%%%%%%%%%%%%%%%%%%%%%%%%%%%%%%%%%%%%%%%%%%%%%%%%%%%%%%%%%%%%%%%%%%%%%%%%

\bibliographystyle{unsrt}
\bibliography{Ref}

\vspace{12pt}
\color{red}

\end{document}